\documentclass[12pt]{article}
\usepackage{mitpress}

\newtheorem{theorem}{Theorem}

\newtheorem{corollary}[theorem]{Corollary}

\newtheorem{lemma}[theorem]{Lemma}

\newtheorem{remark}[theorem]{Remark}

\newenvironment{proof}[1][Proof]{\textbf{#1.} }{\ \rule{0.5em}{0.5em}}
\newdimen\dummy
\dummy=\oddsidemargin
\begin{document}

\title{The fundamental group of the harmonic archipelago}
\author{Paul Fabel \\
%EndAName
Department of Mathematics \& Statistics\\
Mississippi State University}
\maketitle

\begin{abstract}
The harmonic archipelago $HA$ is obtained by attaching a large pinched
annulus to every pair of consecutive loops of the Hawaiian earring. We
clarify $\pi _{1}(HA)$ as a quotient of the Hawaiian earring group, provide
a precise description of the kernel, show that both $\pi _{1}(HA)$ and the
kernel are uncountable, and that $\pi _{1}(HA)$ has the indiscrete topology.
\end{abstract}

\section{Introduction}

This note serves to clarify certain properties of $\pi _{1}(HA)$ and $\pi
_{1}(HE),$ the topological fundamental groups respectively of the harmonic
archipelago and the Hawaiian earring.

The Hawaiian earring $HE\;$is the union of a null sequence of simple closed
curves meeting in a common point.

Introduced by Bogley and Sieradski, the harmonic archipelago $HA,$ is the
space obtained by attaching large pinched annuli, one for each pair of
consecutive loops, to the Hawaiian earring $HE$. In \cite{Bog2} Bogley and
Sieradski construct a comprehensive theory which provides a useful framework
for investigating the fundamental groups of locally complicated spaces such
as $HE$ and $HA.$ Various properties of $\pi _{1}(HA),$ such as its
uncountability, are uncovered in \cite{Bog2}.

In \cite{Biss}, Biss also uses $HE$ and $HA$ as motivating examples, and
proves some nice general results on topological fundamental groups and their
generalized covering spaces.

However, an oversight in \cite{Biss} $($see Remark \ref{counter}), leads to
a false description of $\pi _{1}(HA)$ and its false generalization Theorem
8.1. There is also a typographical error in the description of $\ker
(j^{\ast }),$ the kernel of the epimorphism $j^{\ast }:\pi
_{1}(HE)\rightarrow \pi _{1}(HA),$ induced by inclusion $j:HE\hookrightarrow
HA.$ After adjusting for this, there remains arguably room for further
discussion regarding which elements of $\pi _{1}(HE)$ belong to $\ker
(j^{\ast })$.

For example it follows from the investigations of Morgan/Morrison \cite{mor}%
, DeSmit \cite{DeSmit}, and Cannon/Conner \cite{Cannon} that elements of $%
\pi _{1}(HE)$ can be seen as ``transfinite words'' over an infinite alphabet 
$\{x_{1},x_{2},...\}$ with each letter appearing finitely many times.

As described in \cite{Biss}, $\ker (j^{\ast })$ is generated by the
relations $x_{i}=x_{j}$ for all $i$ and $j.$ Does this mean two transfinite
words over $\{x_{1},x_{2},...\}$ are equivalent in $\pi _{1}(HA)$ if and
only if one can be transformed into the other by finitely many substitutions
and finitely many cancellations of consecutive letters? No, for example $%
x_{1}x_{2}^{-1}(x_{3}x_{4}...)x_{1}x_{2}^{-1}(...x_{4}^{-1}x_{3}^{-1})$ is
trivial in $\pi _{1}(HA),$ but cannot be transformed into the trivial word
with finitely many such operations.

Are two transfinite words over $\{x_{1},x_{2},..\}$ equivalent in $\pi
_{1}(HA)$ if one can be transformed into the other after ``infinitely many
substitutions''? No, for then the essential element $%
(x_{1}x_{2}^{-1}x_{3}x_{4}^{-1}...)$ could be transformed into the
inessential $x_{1}x_{1}^{-1}x_{2}x_{2}^{-1}...$

We provide a precise description of $\ker (j^{\ast })$ in Theorem \ref{main}
and prove as Corollary \ref{cor1} that $\ker (j^{\ast })$ is uncountable.
Corollaries \ref{cor2} and \ref{cor3} provide proofs of results also
indicated in \cite{Bog2} and \cite{Biss}: $\pi _{1}(HA)$ is uncountable and,
despite its large cardinality, $\pi _{1}(HA)$ has the indiscrete topology.

\section{Definitions}

For $n\in \{1,2,3,...\}$ let $X_{n}\subset R^{2}$ denote the circle of
radius $\frac{1}{n}$ centered at $(\frac{1}{n},0).$ Let $Y_{n}=\cup
_{i=n}^{\infty }X_{i}$. Thus $Y_{n}$ is the Hawaiian earring determined by
the loops $X_{n},X_{n+1},...$

Let $A_{n}\subset R^{2}$ denote the closed pinched annulus bounded by $%
X_{n}\cup X_{n+1}.$ Let $Y^{n}=Y_{1}\cup A_{1}..\cup A_{n}.$ Endow $Y^{n}$
with the subspace topology inherited from $R^{2}.$

For the \textbf{underlying set} let $HA=\cup _{n=1}^{\infty }Y^{n}.$ However
we define the topology of $HA$ such that $Y^{n}$ inherits the usual topology
but such that 1) There exists a sequence $z_{n}\in int(A_{n})$ such that $%
\{z_{1},z_{2}...\}$ has no subsequential limit and 2) If $p\in HA$ is a
subsequential limit of the sequence $y_{1},y_{2},..$, $\ $and if $y_{n}\in
int(A_{n})$ for all $n,$ then $p=(0,0)$.

Let $G_{n}=\pi _{1}(Y_{n},(0,0)).$

Let $F_{N,n}$ denote the free group on the letters $\{x_{N},x_{N+1},...x_{n}%
\}$ coupled with the symbol $1$ denoting the trivial element.

Let $\phi _{N,n}:F_{N,n+1}\rightarrow F_{N,n}$ denote the homomorphism such
that $\phi _{N,n}(x_{i})=x_{i}$ if $N\leq i\leq n$ and $\phi
_{n}(x_{n+1})=1. $

For $N\geq 1$ let $G^{N}$ denote the inverse limit of free groups determined
by $F_{N,N}\leftarrow F_{N,N+1}.....$ under the bonding maps $\phi _{N,n}.$

Formally elements of $F_{N,n}$ are equivalence classes of words under the
obvious cancellations, and the group operation is catcatanation. However
each element of $F_{N,n}$ has a unique representative with a minimal number
of nontrivial letters. Consequently each element of $G^{N}$ is uniquely
determined by a \textbf{canonical sequence} $w_{N},w_{N+1},...$ such that $%
w_{n}\in F_{N,n}$ and $w_{n}$ is a maximally reduced word in $F_{N,n}.$

\section{$\protect\pi _{1}(HA)$}

It is a nontrivial\textbf{\ }fact (\cite{Cannon}, \cite{DeSmit} \cite{mor})
that $\pi _{1}(Y_{1})$ injects naturally into the inverse limit of free
groups.

\begin{remark}
\label{morg}Given $n\geq N\geq 1$ let $Z_{N,n}=\cup _{i=N}^{n}X_{n}.$ Let $%
r_{N,n}:Y_{N}\rightarrow Z_{N,n}$ denote the retraction collapsing $X_{i}$
to $(0,0)$ for $i>n.$ Since $r_{N,n}(r_{N,n+1})=r_{N,n},$ the maps $r_{N,n}$
induce a homomorphism $\psi _{N}:\pi _{1}(Y_{N},(0,0))\rightarrow
\lim_{\leftarrow }\pi _{1}(Z_{N,n},(0,0)).$ It is shown in \cite{mor} and 
\cite{DeSmit} that $\psi _{N}$ is one to one. Moreover $G_{N}=im(\psi _{N})$
consists of all canonical sequences $(w_{N},w_{N+1},...)$ such that for each 
$i$ there exists $M_{i}$ such that for all $n\geq N,$ $x_{i}$ appears at
most $M_{i}$ times in $w_{n}.$
\end{remark}

It is falsely asserted in \cite{Biss} that $\psi _{1}$:$\pi
_{1}(HE)\rightarrow \lim_{\leftarrow }F_{1,n}$ is an isomorphism.
Consequently Theorem 8.1 of \cite{Biss} is false.

\begin{remark}
\label{counter}The canonical monomorphism $\psi _{1}:\pi
_{1}(Y_{1},(0,0))\hookrightarrow G^{1}$ is \textbf{not} surjective. By
compactness of $[0,1],$ a given path in $Y_{1}$ can traverse each loop only
finitely many times. Thus the element 
\[
(1,x_{1}x_{2}x_{1}^{-1}x_{2}^{-1},x_{1}x_{2}x_{1}^{-1}x_{2}^{-1}x_{1}x_{3}x_{1}^{-1}x_{3}^{-1},...)\in \lim_{\leftarrow }F_{1,n}
\]
has no preimage in the Hawaiian earring group.
\end{remark}

\begin{remark}
\label{rem1}$X_{n}\cup X_{n+1}$ is a strong deformation retract of $%
A_{n}\backslash \{z_{n}\}.$ Thus $Y^{n}$ is a strong deformation retract of $%
HA\backslash \{z_{n+1},z_{n+2},...\}.$
\end{remark}

\begin{lemma}
\label{lem1}Suppose $f:S^{1}\rightarrow Y_{1}$ is any map. Then $f$ is
inessential in $HA$ if and only if there exists $N$ such that $f$ is
inessential in $Y^{N}.$
\end{lemma}

\begin{proof}
Suppose $f$ is inessential in $HA.$ Let $F:D^{2}\rightarrow HA$ be a
continuous extension of $f.$ Since $im(F)$ is compact, there exists $N$ such
that $z_{n}\notin im(F)$ whenever $n>N.$ Thus $im(F)\subset HA\backslash
\{z_{N+1},z_{N+2},..\}.$ By remark \ref{rem1} $Y^{N}$ is a strong
deformation retract of $HA\backslash \{Z_{N+1},Z_{N+2},...\}.$ Thus $f$ is
inessential in $Y^{N}.$ Conversely if there exists $N$ such that $f$ is
inessential in $Y^{N}$ then $f$ is inessential in $HA$ since $Y^{N}\subset
HA.$
\end{proof}

Notice $Y_{N}$ is a strong deformation retract of $Y^{N-1}$ under a homotopy 
$R_{N,t}:Y^{N-1}\rightarrow Y_{N}$ collapsing $A_{1}\cup ..A_{N-1}$onto the
simple closed curve $X_{N}.$ In particular given any loop $%
f:S^{1}\rightarrow Y_{1}$ we may canonically deform $f$ in $Y^{N-1}$ to a
loop $g$ such that $im(g)\subset Y_{N}.$ Appealing to Remark \ref{morg} we
may identify $\pi _{1}(Y_{i},(0,0))$ with $G_{i}.$ Thus the composition $%
Y_{1}\hookrightarrow Y^{N-1}\rightarrow Y_{N}$ $\ $induces a homomorphism $%
q_{N}^{\ast }:G_{1}\rightarrow G_{N}.$ Combining these observations we
obtain the following:

\begin{lemma}
\label{lem2}$q_{N}^{\ast }(w_{1},w_{2},...)=(v_{N},v_{N+1},..)$ if and only
the following property is satisfied for each $n\geq N$: For each $i\leq N$
replace each occurrence of $x_{i}$ in $w_{n}$ with $x_{N},$ creating a word $%
w_{n}^{\symbol{94}}$ on the letters $\{x_{N},..x_{n}\}.$ Then the word $%
w_{n}^{\symbol{94}}$ is equivalent to $v_{n}$ in the free group $F_{N,n}.$
\end{lemma}

Lemma \ref{onto} is also handled in \cite{Biss} and \cite{Bog2}. Our proof
is similar to an argument that the $2$ sphere is simply connected.

\begin{lemma}
\label{onto}Let $j:Y_{1}\hookrightarrow HA$ denote the inclusion map. Then
the induced homomorphism $j^{\ast }:\pi _{1}(Y_{1},(0,0))\rightarrow \pi
_{1}(HA,(0,0))$ is surjective.
\end{lemma}

\begin{proof}
Suppose $f:S^{1}\rightarrow HA$ is any map such that $f(1)=(0,0)$. Let $J$
be a (nonempty) component of $f^{-1}(HA\backslash Y_{1}).$ Note $%
HA\backslash Y_{1}$ is the union of pairwise disjoint connected open sets $%
int(A_{1})\cup int(A_{2}),...$Thus there exists $i$ such that $f(J)\subset
int(A_{i})$ and $f(\partial J)\subset \partial A_{i}.$ If $z_{i}\in im(f)$
replace $f_{\overline{J}}$ by a path homotopic small perturbation $f_{%
\overline{J}}^{\symbol{94}}$ such that $z_{i}\notin im(f_{\overline{J}}^{%
\symbol{94}})$. By uniform continuity of $f$ finitely many such surgeries
are required. Note $f$ and $f^{\symbol{94}}$ are homotopic in $HA$ and $%
im(f^{\symbol{94}})\subset HA\backslash \{z_{1},z_{2},...\}.$

By remark \ref{rem1} $Y_{1}$ is a strong deformation retract of $%
HA\backslash \{z_{1},z_{2},..\}$ and hence there exists $f^{\ast
}:S^{1}\rightarrow Y_{1}$ such that $j^{\ast }[f^{\ast }]=[f^{\symbol{94}%
}]=[f].$
\end{proof}

Since $j^{\ast }$ is a surjection $\pi _{1}(HA,(0,0))$ is isomorphic to the
quotient group $\pi _{1}(Y_{1},(0,0))/\ker j^{\ast }.$ Combining Lemmas \ref
{lem1} and \ref{lem2} we obtain the following characterization of $\ker
j^{\ast }.$

\begin{theorem}
\label{main}Let $j^{\ast }:G_{1}\rightarrow \pi _{1}(HA,(0,0))$ denote the
epimorphism induced by inclusion $j:Y_{1}\hookrightarrow HA.$ Let $K=\ker
j^{\ast }.$ Then $(w_{1},w_{2},..)\in K$ if and only if there exists $N$
such that the following holds: Suppose for each $i\leq N$ and for each $%
n\geq N$ each occurrence of $x_{i}$ in $w_{n}$ is replaced by $x_{N}$
creating a (nonreduced) word $v_{n}$ on the letters $%
\{x_{N},x_{N+1},...x_{n}\}.$ Then $v_{n}$ can be reduced to the trivial
element of $F_{N,n}.$
\end{theorem}

\begin{corollary}
\label{cor1}$\ker (j^{\ast })$ is uncountable.
\end{corollary}

\begin{proof}
There exist uncountably many distinct permutations of the set $%
\{3,4,5,6,...\}.$ Moreover, the loops in $Y_{3}$ determined by distinct
permutations of $x_{3},x_{4},...$ determine distinct elements of $\pi
_{1}(Y_{3},(0,0))$. Thus $G_{3}$ is uncountable. Each $w\in G_{3}$
determines a homotopically distinct loop in $Y_{1}$ corresponding to $%
x_{1}x_{2}^{-1}wx_{1}x_{2}^{-1}w^{-1}.$ By Theorem \ref{main} this element
is in $\ker (j^{\ast }).$
\end{proof}

Corollary \ref{cor2} is also treated in \cite{Bog2}.

\begin{corollary}
\label{cor2}$\pi _{1}(HA,(0,0))$ is uncountable.
\end{corollary}

\begin{proof}
Consider the set $A$ of functions from $\{1,3,5,...\}\rightarrow \{0,1\}.$
Each element $f\in A$ determines a permutation $\tau
_{f}:\{1,2,3,...\}\rightarrow \{1,2,3,..\}$ as follows: If $f(2n-1)=0$ then $%
\tau _{f}$ fixes $2n-1$ and $2n.$ Otherwise $\tau _{f}$ swaps $2n-1$ and $%
2n. $ The permutation $\tau _{f}$ determines a loop in $Y_{1}$ by the
following recipe. Travel clockwise once around $X_{\tau _{f}(1)},$ then
travel clockwise once around $X_{\tau _{f}(2)},$etc... This corresponds to
the transfinite word $\tau _{f}(1)\tau _{f}(2)...\in \pi _{1}(Y_{1}).$

Note $A$ is uncountable. Declare two elements $\{f,g\}\subset A$ equivalent
if $f$ and $g$ agree except on a finite set. Each equivalence class in $A$
has countably many elements and hence $B,$ the set of equivalence classes of 
$A,$ is uncountable. If $[f]$ and $[g]$ are distinct elements of $B$ then by
Lemma \ref{lem2}, for each $N,$ $\tau _{f}$ and $\tau _{g}$ fail to be
equivalent in $Y^{N}.$ Thus by Lemma \ref{lem1} $\tau _{f}$ and $\tau _{g}$
determine distinct elements of $\pi _{1}(HA,(0,0)).$Thus, since $B$ is
uncountable, $\pi _{1}(HA,(0,0))$ is uncountable.
\end{proof}

It shown in \cite{Biss} (and generalized in \cite{Fabel}) that the quotient
space consisting of path components of based loops in a topological space $X$
is a topological group. Despite a large number of elements, $\pi
_{1}(HA,(0,0))$ has the indiscrete topology. Corollary \ref{cor3} is also
argued in \cite{Biss}.

\begin{corollary}
\label{cor3}The topological group $\pi _{1}(HA,(0,0))$ has exactly one
nonempty open subset.
\end{corollary}

\begin{proof}
Let $L=\{f:[0,1]\rightarrow HA|f(0)=f(1)=(0,0)\}$ with the uniform topology$%
. $ By definition $\pi _{1}(HA,(0,0))$ is the collection of path components
of $L$ with the quotient topology. Suppose $B\subset \pi _{1}(HA(0,0)),$
suppose $B\neq \emptyset $ and suppose $[f]\in B.$

Since $im(f)$ is compact and $\{z_{1},z_{2},..\}$ is not compact, there
exists $N$ such that $im(f)\subset HA\backslash \{z_{N+1},z_{N+2},..\}.$ For 
$M>N$ there is a strong deformation retraction from $HA\backslash
\{z_{M+1},z_{M+2},..\}$ onto $Y^{M}$. Thus for large $M$ we may deform $f$
in $HA$ to a nearby map $f_{M}$ such that $im(f_{M})\subset Y^{M}.$ Let $%
A_{M}=f_{M}^{-1}(Y^{M}\backslash (X_{M+1}\cup X_{M+2}...)).$ Redefine $f_{M}$
over $A_{M}$ to be the constant $(0,0)$ creating a map $f_{M}^{\symbol{94}}.$
Note $im(f_{M}^{\symbol{94}})\subset A_{1}\cup ...A_{M-1}.$ Notice $%
A_{1}\cup ..A_{M-1}$ has the homotopy type of $S^{1}.$ By wrapping around $%
X_{M}$ as many times as necessary, we may extend $f_{M}^{\symbol{94}}$ to a
map of $f_{M}^{\ast }$:$[0,1+\frac{1}{M}]\rightarrow HA$ such that $%
f_{M}^{\ast }$ is an inessential loop and $f(t)\in X_{M}$ for $t\in \lbrack
1,\frac{1}{M}]$. Reparameterize $f_{M}^{\ast }$ linearly to create a map $%
f_{M}^{\ast \ast }:[0,1]\rightarrow HA.$ Notice $f_{M}^{\ast \ast
}\rightarrow f$ uniformly. Thus the path component of the trivial loop is
dense in the space of based loops of $HA.$ Thus, if $B$ is open then $B$
contains the trivial element.

Suppose on the other hand that $B$ is closed. Choose $N$ such that $%
im(f)\subset Y^{n}$ whenever $n\geq N.$ Notice $Y_{n+1}$ is a strong
deformation retract of $Y^{n}.$ Thus $[f]$ contains representatives whose
images have arbitrarily small diameter. Thus $B$ contains the trivial
element.

Thus the trivial element of $\pi _{1}(HA(0,0))$ belongs to each nonempty
open set and each nonempty closed set. Hence $\pi _{1}(HA(0,0))$ has only
one nonempty open subset.
\end{proof}


\begin{thebibliography}{9}
\bibitem{Biss}  Biss, Daniel K. The topological fundamental group and
generalized covering spaces. Topology Appl. 124 (2002), no. 3, 355--371.

\bibitem{Bog2}  Bogley, W.A., Sieradski, A.J. Universal path spaces.
Preprint. http://oregonstate.edu/\symbol{126}bogleyw/\#research.

\bibitem{Cannon}  Cannon, J. W.; Conner, G. R. The combinatorial structure
of the Hawaiian earring group. Topology Appl. 106 (2000), no. 3, 225--271.

\bibitem{DeSmit}  de Smit, Bart. The fundamental group of the Hawaiian
earring is not free. Internat. J. Algebra Comput. 2 (1992), no. 1, 33--37.

\bibitem{Fabel}  Fabel, Paul. Completing Artin's braid group on infinitely
many strands. http://front.math.ucdavis.edu/math.GT/0201303. To appear in
Journal of Knot Theory and its Ramifications.

\bibitem{mor}  Morgan, John W.; Morrison, Ian. A van Kampen theorem for weak
joins. Proc. London Math. Soc. (3) 53 (1986), no. 3, 562--576.
\end{thebibliography}
\end{document}